\documentclass[12pt,a4paper]{article}
\usepackage{amsmath,amsfonts,amssymb,amscd}

\usepackage[a4paper]{geometry}

\def\Ker{\operatorname{Ker}}

\def\id{\operatorname{id}}

\def\Im{\operatorname{Im}}

\newcounter{th}
\def\t{\refstepcounter{th}{\bf \noindent{Theorem} \arabic{th}. }}

\newcounter{prop}
\def\prop{\refstepcounter{prop}{\bf \noindent{Proposition} \arabic{prop}. }}

\newcounter{lem}

\newcounter{de}

\newcounter{ex}

\begin{document}

\begin{center}
    {\LARGE{\bf Geometric approach to  $p$-singular Gelfand-Tsetlin $\mathfrak {gl}_n$-modules }}
\end{center}
\medskip

\begin{center}
{\large Elizaveta Vishnyakova}
\end{center}

\bigskip

\bigskip

\begin{abstract}	
We give an elementary construction of a $p\geq 1$-singular Gelfand-Tsetlin $\mathfrak{gl}_n(\mathbb C)$-module in terms of local distributions. This is a generalization of the universal $1$-singular Gelfand-Tsetlin $\mathfrak{gl}_n(\mathbb C)$-module obtained in \cite{Futorny}. We expect that the family of new Gelfand-Tsetlin modules that we obtained will lead to a classification of all irreducible $p>1$-singular Gelfand-Tsetlin modules. So far such a classification is known only for singularity $n=1$.

\end{abstract}

\section{Introduction}
In classical Gelfand-Tsetlin Theory one constructs explicitly an action of the Lie algebra $\mathfrak g:=\mathfrak{gl}_n(\mathbb C)$ in a basis forming by Gelfand-Tsetlin tableaux. Let $V$ be the vector space of all Gelfand-Tsetlin tableaux of fixed order, see the main text for details, and $\mathcal U(\mathfrak g)$ be the universal enveloping algebra of $\mathfrak g$. In \cite{Fu-Ovs} it was proved that there exists a ring structure on the vector space $\mathcal R:=H^0(V,\mathcal M\star \Gamma)$, where $\mathcal M\star \Gamma$ is the sheaf of meromorphic functions on $V$ with values in a certain group $\Gamma$ acting on $V$,  such that  the classical Gelfand-Tsetlin formulas define a homomorphism of rings $\Phi:	\mathcal U(\mathfrak g) \to \mathcal R$.  In the case when $\Im\Phi$  is holomorphic on a certain orbit $\Gamma(o)$ of a point $o\in V$  the	converse statement is also true: {\it any homomorphism of rings $\Phi:	\mathcal U(\mathfrak g) \to \mathcal R$ defines a Gelfand-Tsetlin like formulas with basis forming by elements of $\Gamma(o)$.}  The study of the case when $\Im\Phi$  is not holomorphic in $\Gamma(o)$, but one singular, was initiated  in \cite{Futorny}. For instance the authors \cite{Futorny} constructed the universal $1$-singular Gelfand-Tsetlin $\mathfrak{gl}_n(\mathbb C)$-module using pure algebraic methods. Another version of the construction from \cite{Futorny} of $1$-singular Gelfand-Tsetlin $\mathfrak{gl}_n(\mathbb C)$-module can be found in \cite{Zad}.

In the present paper we study the case of singularity $p\geq 1$. In fact, we give a new elementary geometric construction  of $p$-singular Gelfand-Tsetlin $\mathfrak{gl}_n(\mathbb C)$-modules. In the case $p=1$ our construction gives another version of the constructions of the universal $1$-singular Gelfand-Tsetlin $\mathfrak{gl}_n(\mathbb C)$-modules obtained  in \cite{Futorny} and \cite{Zad}. 
The universal $1$-singular Gelfand-Tsetlin $\mathfrak{gl}_n(\mathbb C)$-modules was used in \cite{Futorny,Futorny2} to classify all irreducible Gelfand-Tsetlin modules with $1$-singularity. We expect that our construction of $p$-singular Gelfand-Tsetlin $\mathfrak{gl}_n(\mathbb C)$-modules will lead to a classification of all irreducible $p>1$-singular Gelfand-Tsetlin modules.

Moreover, our approach leads to a geometric expla\-nation of the formulas from \cite{Futorny,Zad} for the $\mathfrak{gl}_n(\mathbb C)$-action in the universal $1$-singular Gelfand-Tsetlin basis. 
On the other side, our approach may be used for other homomorphisms $\Psi:	\mathcal U(\mathfrak h) \to \mathcal R$, where $\mathfrak h$ is any Lie algebra. In fact, we describe necessary conditions for the image of the homomorphism $\Psi:	\mathcal U(\mathfrak h) \to \mathcal R$ such that certain local distributions supported on the elements of the orbit $\Gamma(o)$ form a basis for an $\mathfrak h$-module.

\section{Preliminaries}

 Let $V$ be the vector space of Gelfand-Tsetlin tableaux $V\simeq \mathbb C^{n(n+1)/2}=\{ (x_{ki})\,\, |\,\,  1\leq i \leq k \leq n  \}$, where $n\geq 2$, and $\Gamma \simeq \mathbb C^{n(n-1)/2}$ be the free abelian group generated by $\sigma_{st}$, where $1\leq t \leq s \leq n-1$.  We fix the following action of $\Gamma$ on $V$: $\sigma_{st} (x) = (x_{ki} +\delta_{ki}^{st})$, where $x = (x_{ki})\in V$ and $\delta_{ki}^{st}$ is the  Kronecker delta. We put $G=S_1\times S_2\times \cdots \times S_n$, where $S_i$ is the symmetric group of degree $i$. The group $G$ acts on $V$ in the following way $(s (x))_{ki} = x_{k s_k(i)}$, where $s=(s_1,\ldots, s_n)\in G$. 

Denote by $\mathcal M$ and by $\mathcal O$ the sheaves of meromorphic and holomorphic functions on $V$, respectively. Let us take $f\in H^0(V,\mathcal M)$, $s\in G$ and $\sigma\in \Gamma$. We set 
$$
s(f) = f\circ s^{-1},\quad  \sigma(f) = f\circ \sigma^{-1} \quad \text{and} \quad s(\sigma) = s\circ \sigma\circ s^{-1}.
$$ 
 Denote by $\mathcal M\star\Gamma: = \bigoplus\limits_{\sigma\in \Gamma} \mathcal M\sigma $ the sheaf of meromorphic functions on $V$ with values in $\Gamma$. In other words, $\mathcal M\star \Gamma$ is the sheaf of meromorphic sections of the trivial bundle $V\times \bigoplus\limits_{\sigma\in \Gamma} \mathbb C\sigma \to V$. An element in $\mathcal M\star\Gamma$ is a finite linear combination of $f\sigma$, where $f\in \mathcal M$ and $\sigma\in \Gamma$. There exists a natural structure of a skew group ring on $H^0(V,\mathcal M\star \Gamma)$, see \cite{Fu-Ovs}. Indeed, 
$$
\sum_if_i \sigma_i \circ 
\sum_j f'_j \sigma'_j := 
\sum_{ij} f_i \sigma_i(f'_j) \sigma_i\circ \sigma'_j.
$$
Here $f_i, f'_j\in H^0(V,\mathcal M)$  and $\sigma_i,\sigma'_j \in \Gamma$. This skew ring we will denote by $\mathcal R$. To simplify notations we use $\circ$ for the multiplication in $\mathcal R$ and for the product in $\Gamma$. 
 We will also consider the multiplication $A*B := B\circ A$ in $H^0(V,\mathcal M\star \Gamma)$. 
The ring $\mathcal R$ possesses the following action of the group $G$
$$
s\big(\sum_if_i \sigma_i\big) = \sum_i s( f_i) s(\sigma_i). 
$$
It is easy to see that this action preserves the multiplication in $\mathcal R$. Hence the vector space $\mathcal R^{G}$ of all $G$-invariant elements is a subring in $\mathcal R$.

The classical Gelfand-Tsetlin formulas  have the following form in terms of generators:
\begin{equation}\label{eq G-Ts generators}
\begin{split}
E_{k,k+1}(T(v)) &= - \sum_{i=1}^k \frac{\prod_{j= 1}^{k+1} (x_{ki}- x_{k+1,j})}{\prod_{j\ne i}^k (x_{ki}- x_{kj})} T(v + \delta_{ki});\\
E_{k+1,k} (T(v)) & = \sum_{i=1}^k \frac{\prod_{j= 1}^{k-1} (x_{ki}- x_{k-1,j})}{\prod_{j\ne i}^k (x_{ki}- x_{kj})} T(v - \delta_{ki});\\
E_{k,k} (T(v)) & = \Big( \sum_{i=1}^k (x_{ki} +i-1) - \sum_{i=1}^{k-1} (x_{k-1,i} +i-1)  \Big) (T(v)),
\end{split}
\end{equation}
see for instance \cite{Futorny}, Theorem $3.6$. Here $E_{st}\in \mathfrak{gl}_n(\mathbb C)$, $T(v)\in V$ is a point in $V$ with coordinates $v=(x_{ki})$ and $T(v \pm \delta_{ki})\in V$ is the tableau obtained by adding $\pm 1$ to the $(k,i)$-th entry of $T(v)$. A Gelfand-Tsetlin tableau is called {\it generic} if $x_{rt} - x_{rs} \notin \mathbb Z$ for any $r$ and for any $s\ne t$. In the case when  $T(v)$ is a generic Gelfand-Tsetlin tableau Formulas (\ref{eq G-Ts generators}) define a $\mathfrak{gl}_n(\mathbb C)$-module structure on the vector space spanned by the elements of the orbit $\Gamma(T(v))$, see for instance Theorem $3.8$ in \cite{Futorny} and references therein. Note that the action of $\Gamma$ in $V$ is free. Hence, $\Gamma(T(v))\simeq \Gamma$ and the elements of the orbit $\Gamma(T(v))$ form the Gelfand-Tsetlin basis. Another observation is that the coefficients in Formulas (\ref{eq G-Ts generators}) are holomorphic in sufficiently small neighborhood of $\Gamma(T(v))$ for a generic $T(v)$.

 There is a natural action of $\mathcal R$ on $H^0(V,\mathcal M)$ that is given by 
 $$
 F \mapsto (f\sigma)(F):= f \sigma (F) = f F \circ \sigma^{-1}.
 $$ 
  Let us identify $T(v)\in V$ with the corresponding evaluation map $ev_v: H^0(V,\mathcal O) \to \mathbb C$, $F\mapsto F(v)$. (Note that $V$ is a Stein manifold, so this identification exists.) Assume that $R = f_i\sigma_i \in \mathcal R$ is holomorphic in a neighborhood of $v$. Then we have $ev_v\circ f_i\sigma_i(F) = f_i(v) F(\sigma_i^{-1} (v))$, where $F\in H^0(V,\mathcal O)$. We put $R(ev_v) := ev_v\circ R$. In these notations we have $R_1(R_2(ev_v)) = ev_v \circ (R_2\circ R_1) = (R_1* R_2)(ev_v)$, where $R_1,R_2 \in \mathcal R$. Now we can rewrite Formulas (\ref{eq G-Ts generators}) in the following form 
 \begin{equation}\label{eq classical formulas G-Ts}
 \Phi (E_{st}) (ev_v) = ev_v \circ \Phi(E_{st}).
 \end{equation}
  Here $\Phi(E_{st})\in \mathcal R$ is defined by Formulas (\ref{eq G-Ts generators}). For example,
 $$
 \Phi (E_{k,k+1}) =  - \sum_{i=1}^k \frac{\prod_{j= 1}^{k+1} (x_{ki}- x_{k+1,j})}{\prod_{j\ne i}^k (x_{ki}- x_{kj})} \sigma^{-1}_{ki}.
 $$ 
In this formula we interpret $x_{ki}$ as the coordinate functions on $V$. Now we see that the statement of Theorem $3.8$ in \cite{Futorny} is equivalent to
 \begin{equation}\label{eq Theorem 38 is equivalent}
(\Phi(X)* \Phi(Y))(ev_v) - (\Phi(Y)* \Phi (X))(ev_v) = \Phi( [X,Y] )(ev_v),
 \end{equation}
 where $X,Y\in \mathfrak g$ and $v$ is generic.  In \cite{Fu-Ovs} the following theorem was proved.

 \medskip
 
 \t\label{teor Fut Ovs}[Futorny-Ovsienko] {\sl The classical Gelfand-Tsetlin formulas (\ref{eq G-Ts generators}) define a homomorphism of rings $\Phi:	\mathcal U(\mathfrak g_{-}) \to \mathcal R$, where $\Phi(X)$, $X\in \mathfrak g$, is as in (\ref{eq classical formulas G-Ts}) and $\mathfrak g_{-} = \mathfrak {gl}_n(\mathbb C)$ with the multiplication $[X,Y]_{-} = -[X,Y] = Y\circ X - X\circ Y$. 
 	
 }
 
 \medskip
 
 \noindent{\it Proof.} Let us give a proof of this theorem for completeness using complex analysis. First of all define a homomorphism of the (free associative) tensor algebra $\mathcal T(\mathfrak g)$ to $\mathcal R$ using (\ref{eq G-Ts generators}). Such a homomorphism always exists because $\mathcal T(\mathfrak g)$ is free. We need to show that the ideal generated by the relation $\Phi(Y) \circ \Phi(X) - \Phi (X) \circ \Phi(Y) -  \Phi([X,Y])$ maps to $0$ for any $X,Y\in \mathfrak g$. In fact we can rewrite Formulas (\ref{eq Theorem 38 is equivalent}) in the following form for any generic $v\in V$, any $X,Y\in \mathfrak g$ and any $F\in H^0(V,\mathcal O)$:
 \begin{align*}
 ev_v \circ (\Phi(Y) \circ \Phi(X)(F) - \Phi(X) \circ \Phi(Y)(T)) = ev_v \circ \Phi([X,Y]) (F).
  \end{align*}
 Since generic points $v$ are dense in $V$, the following holds for any holomorphic $F$:
 \begin{align*}
(\Phi(Y) \circ \Phi(X)(F) - \Phi(X) \circ \Phi(Y)(F)) =  \Phi([X,Y]) (F).
 \end{align*}
 It is remaining to prove that if $R\in \mathcal R$ such that $R(F)=0$ for any $F\in H^0(V,\mathcal O)$, then $R=0$. Indeed, let $R=\sum\limits_{i=1}^s f_i\sigma_i$ and $U\subset V$ be a sufficiently small open set. Clearly it is enough to prove a local version of our statement for any such $U$: from $R(F)|_U=0$ for any $F\in H^0(V,\mathcal O)$ and $R\in \mathcal R|_{U}$, it follows that $R|_{U}=0$.

 Firstly assume that all functions $f_i$ are holomorphic in $U$ and $x_0\in U$.  Let us fix $i_0\in \{1,\ldots,s\}$ and let us take  $F\in H^0(V,\mathcal O)$ such that $F(\sigma_{i_0}^{-1}(x_0))\ne 0$ and $F(\sigma_i^{-1}(x_0))= 0$ for $i\ne i_0$. Then $R(F)(x_0) = (\sum\limits_{i=1}^s f_i\sigma_i(F)) (x_0) = f_{i_0}(x_0)F(\sigma_{i_0}^{-1}(x_0))=0$. Hence $f_{i_0}(x_0)=0$. Therefore, $f_{i_0}|_{U}=0$ for any $i_0$ and $R|_U=0$. 
  
 Further, by induction assume that our statement holds for $(q-1)$ non-holomor\-phic in $U$ coefficients $f_1, \ldots, f_{q-1}$, where $q-1<s$.  Consider meromorphic in $U$ functions $f_1=g_1/h_1, \ldots, f_q$, where $g_1,h_1$ are holomorphic  in $U$ without common non-invertible factors. Then $h_1 R\in \mathcal R|_{U}$ satisfies the equality $h_1 R (F)|_U=0$ for any $F$ and it has $(q-1)$ non-holomorphic summands. Therefore, $h_1 R|_U=0$ and in particular $g_1=0$.$\Box$ 
 
\medskip
 
 \noindent{\bf Remark.} It is well-known that the image $\Phi(\mathcal U(\mathfrak g_{-}))$ is $G$-invariant. This fact can be also verified directly. 
 
 \medskip

The interpretation of a point $T(v)$ as an evaluation map  $ev_v$ suggests a possibility to define a $\mathfrak{gl}_n(\mathbb C)$-module structure on local distributions, i.e. on linear maps $D_v: \mathcal O_v\to \mathbb C$ with $\mathfrak m_v^s\subset \Ker(D_v)$, where $s>0$ and  $\mathfrak m_v$ is the maximal ideal in the local algebra $\mathcal O_v$. In \cite{Futorny} the authors consider formal limits $\lim\limits_{v\to v_0}(T(v+z) - T(v+\tau(z)))/(x_{ki} - x_{kj})$, where $\tau\in G$ is a certain involution, $i\ne j$, $v_0\in V$ is an $1$-singular tableau, see Section $5$, and $v$ is a generic tableau. In fact this limit may be interpreted as a sum of local distributions, see Section $5$. However geometric interpretations were not given in \cite{Futorny}. 
The idea to use local distributions we develop in the present paper. In more details, let $R\in \mathcal R$ and $D_v$ be a  local distribution. We have an action of $(\mathcal R,*)$ on  local distributions defined by $R(D_v)  =   D_v \circ R$. Indeed,
$$
(R_1* R_2)(D_v) = D_v \circ (R_1* R_2) = D_v \circ R_2\circ   R_1  = R_1(R_2(D_v))  , \quad R_i\in \mathcal R.
$$
 By Theorem \ref{teor Fut Ovs} we have $(\Phi(X)* \Phi(Y))(D_v) - (\Phi(Y)* \Phi (X))(D_v) = \Phi( [X,Y] )(D_v)$, if this expression is defined, 
  and  $D_v \mapsto \Phi(X)(D_v)$ gives a structure of a $\mathfrak{gl}_n(\mathbb C)$-module on local distributions. 
Our goal now is to find orbits $\Gamma(o)$ and describe local distributions at points of $\Gamma(o)$ such that this formula is defined.

\section{Gelfand-Tsetlin modules}

In this section we follow \cite{Futorny,Zad}. Consider the universal enveloping algebra $\mathcal U(\mathfrak g)$ of $\mathfrak g = \mathfrak {gl}_n(\mathbb C)$. We have the following sequence of subalgebras
\begin{align*}
\mathfrak {gl}_1(\mathbb C) \subset \mathfrak {gl}_2(\mathbb C) \subset \cdots \subset  \mathfrak {gl}_n(\mathbb C).
\end{align*}
This sequence induces the sequence of the corresponding enveloping algebras 
$$
\mathcal U(\mathfrak {gl}_1(\mathbb C)) \subset \cdots \subset  \mathcal U(\mathfrak {gl}_n(\mathbb C)). 
$$
Denote by $Z_m$ the center of $\mathcal U(\mathfrak {gl}_m(\mathbb C))$, where $1\leq m \leq n$. The subalgebra $\Upsilon$ in $\mathcal U(\mathfrak {gl}_n(\mathbb C))$ generated by elements of $Z_m$, where $1\leq m \leq n$, is called the {\it Gelfand-Tsetlin subalgebra of $\mathcal U(\mathfrak {gl}_n(\mathbb C))$}. This subalgebra is the polynomial algebra with $n(n+1)/2$ generators $(c_{ij})$, where $1\leq j \leq i \leq n$, see \cite{Futorny}, Section $3$. Explicitly these generators are given by
$$
c_{ij} = \sum\limits_{(s_1,\ldots, s_j)\in \{1,\ldots,i\}^j} E_{s_1s_2} E_{s_2s_3} \cdots E_{s_ks_1},
$$
where $E_{st}$ form the standard basis of  $\mathfrak {gl}_n(\mathbb C)$. 

\medskip

\noindent{\bf Definition.} [Definition $3.1$, \cite{Futorny}]  A finitely generated $\mathcal U(\mathfrak {gl}_n(\mathbb C))$-module $M$ is called a {\it Gelfand-Tsetlin module with respect to $\Upsilon$}  if $M$ splits into a direct sum of $\Upsilon$-submodules:
$$
M= \bigoplus_{\mathfrak m } M(\mathfrak m),
$$
where the sum is taken over all maximal ideals $\mathfrak m$ in $\Upsilon$. Here 
$$
M(\mathfrak m) = \{ v\in M \,\,|\,\, \mathfrak m^q(v) = 0 \,\, \text{for some} \,\, q\geq 0 \}.
$$

For  the following theorem we refer \cite{Futorny}, Section $3$ and \cite{Zad}, Theorem $2$. Compare also with Theorem \ref{teor Fut Ovs}. 

\medskip

\t\label{teor G-T modules} {\sl The image $\Phi(\Upsilon)$  coincides with the subalgebra of polynomials in $H^0(V,\mathcal O^G)$. In other words, for any $X\in \Upsilon$ 
we have $\Phi(X) = F \id $, where $F$ is a $G$-invariant polynomial. 

}

\medskip

\noindent{\bf Corollary.} {\sl All modules corresponding to the homomorphism $\Phi$ with a basis forming by local distributions on $V$ are  Gelfand-Tsetlin modules. 

}

\medskip

\noindent{\it Proof.} Indeed, let $D_v$ be a local distribution on $V$. Them for any $X\in \Upsilon$ we have 
$$
\Phi(X)(D_v) =  D_v \circ \Phi(X) =  D_v \circ (F \id).
$$
By definition of a local distribution, $D_v$ annihilates  $\mathfrak m_v^q$ for some $q >0$, where $\mathfrak m_v$ is the maximal ideal in $H^0(V,\mathcal O)$. In particular, $D_v$ annihilates a degree of the corresponding to $\mathfrak m_v$ maximal ideal in $\Phi(\Upsilon)$.$\Box$ 

\medskip

\section{Alternating holomorphic functions}

 An {\it alternating polynomial} is a polynomial $f(x_1,\ldots, x_n)$ such that 
 $$
 f(\tau(x_1),\ldots,\tau(x_n))= (-1)^{\tau} f(x_1,\ldots, x_n),
 $$
 for any $\tau\in S_n$. An example of an alternating polynomial is the {\it Vandermonde determinant} 
 $$
\mathcal V_n = \prod\limits_{1\leq i<j \leq n} (x_j-x_i).
 $$
 In fact this example is in some sense unique. More precisely, 
we need the following property of alternating polynomials.
 
\medskip

\prop\label{prop alternating polynomial} {\sl 
Any alternating polynomial $f(x_1,\ldots, x_n)$ can be written in the form $f=  \mathcal V_n \cdot g$, where $ g=g(x_1,\ldots, x_n)$ is a symmetric polynomial. 

}
 
\medskip
 
 \noindent{\it Proof.} The proof follows from the following facts. Firstly every alternating polynomial $f$ vanishes on the subvariety $x_i=x_j$, where $i\ne j$. Hence $(x_i-x_j)$ is a factor of $f$ and therefore $\mathcal V_n $ is also a factor of $f$. Secondly it is clear that the ratio $g=f/\mathcal V_n$ is a symmetric polynomial.$\Box$ 
 
\medskip

\noindent {\bf Corollary.} {\sl Let $F= F(x_1,\ldots, x_n)$ be a holomorphic alternating function, i.e $ F(\tau(x_1),\ldots,\tau(x_n))= (-1)^{\tau} F(x_1,\ldots, x_n)$  for any $\tau\in S_n$. Then $F=\mathcal V_n  \cdot G$,  
where $G = G(x_1,\ldots, x_n)$ is a symmetric holomorphic function.$\Box$ 

}

\medskip

We will use this Corollary in Gelfand-Tsetlin Theory. Let $o=(x^0_{ki})\in V$ be a  Gelfand-Tsetlin tableau such that $x^0_{ki_1} = \cdots = x^0_{ki_p}$, where $p \geq 2$. Denote by $W$ a sufficiently small neighborhood of the orbit $\Gamma (o)$.  For simplicity we put $x_j:= x_{ki_j}$. Let $\mathcal V_p = \mathcal V(x_1,\ldots,x_p)$ and $S_p\subset G$ be the permutation group of $(x_1,\ldots,x_p)$. We need the following proposition.

\medskip
\prop\label{prop A_1 ...A_n is 1 sing} { \sl Let $A_j = \sum\limits_i (H^j_{i}/\mathcal V_p)\sigma_{i}\in \mathcal R$, where $j=1,\ldots,q$, and $H^j_{i}$ are holomorphic in $W$, be $S_p$-invariant elements in $\mathcal R$.  Then  $A_1\circ \cdots \circ A_q = \sum\limits_r (G_{r}/\mathcal V_p)\sigma_{r}$, where $G_{r}$ are holomorphic at $o$.

}

\medskip

\noindent{\it Proof.}  Assume by induction that for $k=q-1$ our statement holds. In other words, assume that $A_1\circ \cdots \circ A_{q-1} = \sum\limits_i(G_{i}/\mathcal V_p)\sigma_{i}$, where $G_{i}$ are holomorphic at $o$.  We have
\begin{align*}
A_1\circ \cdots \circ A_{q} = \sum\limits_{i,j} \frac{G_{i} \sigma_i(H^q_{j})}{\mathcal V_p\sigma_i(\mathcal V_p)} \sigma_{i} \circ \sigma_{j}.  
\end{align*}
Assume that $(G_{i_0}\sigma_{i_0}(H^q_{j})) /\sigma_{i_0}(\mathcal V_p) $ is singular at $o$. Note that $G_{i}\sigma_i(H^q_{j})$ is holomorphic at $o$. Let $\sigma_{i_0}(x_1,\ldots, x_p) = (x_1+m_1,\ldots, x_p+m_p)$, where $m_i\in \mathbb Z$. Hence, 
$$
\sigma_{i_0}(\mathcal V_p) = \prod\limits_{1\leq i<j \leq p} (x_j-x_i + m_j-m_i).  
$$
If $(G_{i_0}\sigma_{i_0}(H^q_{j})) /\sigma_{i_0}(\mathcal V_p) $ is singular at $o$, we have $m_{i_1} = \cdots =  m_{i_r}$, where $r\leq p$. For instance, $\tau(\sigma_{i_0}) = \sigma_{i_0}$ for any $\tau \in S_r$, where $S_r$ is the permutation group of $(x_{i_1}, \ldots, x_{i_r})$. 
The product $\sum\limits_i(G_{i}/\mathcal V_p)\sigma_{i}$ is $S_p$-invariant since $A_s$ are $S_p$-invariant by assumption. Since $\sigma_{i_0}$ is $S_r$-invariant and the decomposition $\sum\limits_i(G_{i}/\mathcal V_p)\sigma_{i} \in \mathcal R$ is unique, the function $G_{i_0}/\mathcal V_p$ is $S_r$-invariant. Therefore the holomorphic function $G_{i_0}$ is $S_r$-alternating.  By Corollary of Proposition \ref{prop alternating polynomial}, we have $G_{i_0} = \mathcal V_r(x_{i_1}, \ldots, x_{i_r}) G'_{i_0}$. Hence, $(G_{i_0}\sigma_{i_0}(H^q_{j})) /\sigma_{i_0}(\mathcal V_p) $ is holomorphic at $o$.$\Box$ 

\medskip

\section{Main results}

Let $o=(x^0_{ij})$ be a  Gelfand-Tsetlin tableau such that $x^0_{ki_1} = \cdots = x^0_{ki_p}$, where $p \geq 2$, and such that $x^0_{st}-x^0_{sr}\notin \mathbb Z$ otherwise. We will call the orbits  $\Gamma(o)$ of such points {\it $p$-singular} and corresponding modules {\it $p$-singular Gelfand-Tsetlin $\mathfrak {gl}_n(\mathbb C)$-modules}. 
The stabilizer $G_o\subset G$ of $o$ is isomorphic to the permutation group $G_o\simeq S_p$. We put $z_{rt} = x_{ki_r} - x_{ki_t}$, where $r\ne t$. Then $\frac{\partial}{\partial z_{rt}} = \frac12 (\frac{\partial}{\partial x_{ki_r}} - \frac{\partial}{\partial x_{ki_t}})$ are the corresponding derivations. Let us fix a sufficiently small neighborhood $W$ of the orbit $\Gamma(o)$ such that $W$ is $\Gamma$ and  $G$-invariant. 
We put 
\begin{align*}
	\mathcal L:= ev_o \circ \frac{\partial}{\partial z_{12}}\circ \cdots \circ \frac{\partial}{\partial z_{1p}} \circ \frac{\partial}{\partial z_{23}}  \circ \cdots \circ \frac{\partial}{\partial z_{p-1,p}} \cdot  z_{12} \cdots z_{1p} z_{23} \cdots z_{p-1,p}. 
\end{align*}
Clearly $\mathcal L$ is $S_p$-invariant. Indeed, $z_{12} \cdots z_{1p} z_{23} \cdots z_{p-1,p}$ is equal in fact to the Vandermonde determinant $\mathcal V_p$ in $x_{ki_1}, \ldots,x_{ki_p}$. Hence it is alternating, see Section $4$. On the other side, $\tau(\frac{\partial}{\partial z_{rt}}) = \tau\circ \frac{\partial}{\partial z_{rt}} \circ \tau^{-1} = \frac{\partial}{\partial z_{\tau(r)\tau(t)}} $. Therefore the sequence of derivations in $\mathcal L$ is also alternating. 
Further we put 
$$
T:= \{(12), \ldots (1p), (23),\ldots, (p-1,p) \}
$$
and $I,J$ are two subsets in $T$ such that $I\cup J=T$ and $I\cap J=\emptyset$. Then the elements $z_T$, $z_I$, $z_J$ and the elements $\frac{\partial}{\partial z_{T}}$, $\frac{\partial}{\partial z_{I}}$, $\frac{\partial}{\partial z_{J}}$ are the product and the composition of the corresponding $z_{ij}$, respectively. In this notations $\mathcal L= ev_o \circ \frac{\partial}{\partial z_{T}} \cdot z_{T}$. Note that $z_T=\mathcal V_p$. For any subset $I\subset T$ and $\sigma_i\in \Gamma$, consider the following sum of local distributions
\begin{align*}
	D_{I,\sigma_i} :=  \mathcal L \circ \sum_{\tau\in S_p} (-1)^{\tau}\tau(z_I\sigma_i)/z_T.
\end{align*}
Clearly, we have the following relations 
\begin{equation}\label{eq relation for basis general}
	D_{I,\sigma_i} = (-1)^{\tau'} D_{\tau'(I),\tau'(\sigma_i)}, \quad \tau'\in S_p.
\end{equation}
And we do not have other relations here. We need the following proposition. 

\medskip

\prop\label{prop formulas n=1} {\sl Let us take $\sum\limits_ih_i\sigma_i\in \mathcal R$ an  $S_p$-invariant element that satisfies conditions of Proposition \ref{prop A_1 ...A_n is 1 sing} at $o$. Then we have the following equality of holomorphic operators
	\begin{equation}\label{eq formula n=1, sing}
		\begin{split}
			\mathcal L \circ (\sum\limits_i h_i\sigma_i) = 
		\frac{1}{|S_p|}	\sum_{I,i} \frac{\partial g_i}{\partial z_{I}}(o) D_{I,\sigma_i},
		\end{split}
	\end{equation}
where $g_i= z_Th_i$ and $T$ and $I$ are as above. In other words Formula (\ref{eq formula n=1, sing}) means that $(\sum\limits_i h_i\sigma_i) (\mathcal L)$ is a linear combination of the local distributions $D_{I,\sigma_i}$

}

\medskip
\noindent{\it Proof.}  We have
\begin{align*}
\mathcal L(h_i\sigma_i) = ev_o \circ \frac{\partial}{\partial z_{T}} \circ g_i\sigma_i =  ev_o \circ&\Big(\sum_{I} \frac{\partial g_i}{\partial z_{I}}  \frac{\partial}{\partial z_{J}} \circ \sigma_i \Big) = \\
\sum_{I} \frac{\partial g_i}{\partial z_{I}} (o) \Big(ev_o \circ \frac{\partial}{\partial z_{J}} \frac{\partial}{\partial z_{I}}\cdot z_I \circ \sigma_i \Big) = \sum_{I} \frac{\partial g_i}{\partial z_{I}} (o)& \Big(ev_o \circ \frac{\partial}{\partial z_{T}} \cdot z_T \circ\frac{\sigma_i}{z_J}  \Big) =\\
&\sum_{I} \frac{\partial g_i}{\partial z_{I}} (o) \Big(\mathcal L \circ\frac{\sigma_i}{z_J}  \Big).
\end{align*} 
Since $\mathcal L(h_i\sigma_i)$ is $S_p$-invariant, we have:
\begin{align*}
|S_p| \mathcal L(h_i\sigma_i) =  \sum_{I} \frac{\partial g_i}{\partial z_{I}} (o) \mathcal L \circ \Big(\sum_{\tau\in S_p}\frac{(-1)^{\tau}\tau(z_I\sigma_i)}{z_T}  \Big) =  \sum_{I} \frac{\partial g_i}{\partial z_{I}} (o) \mathcal L \circ  D_{I,\sigma_i}.
\end{align*}
The proof is complete.$\Box$

\medskip

\t\label{teor main n>1 gen}{\bf [Main result 1]} {\sl 
	 Let $\mathfrak g$ be any Lie algebra and $\Phi : \mathcal U(\mathfrak g) \to \mathcal R$ be a homomorphism of rings.
	If $\Phi(\mathfrak g)$ is generated by elements satisfying the conditions of Proposition \ref{prop A_1 ...A_n is 1 sing}, then the vector space spanned by the elements $D_{I,\sigma_i}$, where $I\subset T$ is a subset and $\sigma_i\in \Gamma$, up to relations (\ref{eq relation for basis general})  form a basis for the $\mathfrak g$-module.

}

\medskip

\noindent {\it Proof.} Assume that $\Phi(\mathfrak g)$ is generated by $A_i$ as in  Proposition \ref{prop A_1 ...A_n is 1 sing}. By this proposition we see that any product of such generators has the form   $\sum\limits_i (G_{i}/\mathcal V_p)\sigma_{i}$, where $G_{i}$ are holomorphic at $o$.
 We need to prove that $\sum\limits_i (G_{i}/\mathcal V_p)\sigma_{i} (D_{I,\sigma_j})$ is a linear combination of $D_{I',\sigma_{i'}}$, We have
 \begin{align*}
 \sum\limits_i (G_{i}/\mathcal V_p)\sigma_{i} (D_{I,\sigma_j}) = D_{I,\sigma_j} \circ \sum\limits_i (G_{i}/\mathcal V_p)\sigma_{i} =\\
  \mathcal L \circ \Big(\sum_{\tau\in S_p}\frac{(-1)^{\tau}\tau(z_I\sigma_j)}{\mathcal V_p}  \Big)  \circ \sum\limits_i (G_{i}/\mathcal V_p)\sigma_{i}. 
 \end{align*}
Now we apply Proposition \ref{prop A_1 ...A_n is 1 sing} to the composition in the last line. The result follows from Proposition \ref{prop formulas n=1}.$\Box$

\medskip

Note that the classical Gelfand-Tsetlin generators (\ref{eq G-Ts generators}) satisfy conditions of Proposition \ref{prop A_1 ...A_n is 1 sing}. 

\medskip

\t\label{teor main n>1 class}{\bf [Main result 2]} {\sl The vector space spanned by the elements $D_{I,\sigma_i}$, where $I\subset T$ is a subset and $\sigma_i\in \Gamma$, up to relations (\ref{eq relation for basis general})  form a basis for the $\mathfrak{gl}_n(\mathbb C)$-module. By Corollary of Theorem \ref{teor G-T modules}, this module is a Gelfand-Tsetlin module. 	

}

\medskip
\noindent{\bf Remark.} The basis of elements $D_{I,\sigma_i}$, where $I\subset T$ is a subset and $\sigma_i\in \Gamma$, up to relations (\ref{eq relation for basis general}), can be simplified in some cases. In other words sometimes we can find a natural submodules in the corresponding module. It depends on the singularity type of $\Phi(\mathcal U(\mathfrak g))$ at $o$. It is required in Theorem \ref{teor main n>1 gen} than the singularity type is not more than $\mathcal V_p$. However, if the singularity type of $\Phi(\mathcal U(\mathfrak g))$ at $o$ is less than type of $\mathcal V_p$, Formula (\ref{eq formula n=1, sing})  says  that we can reduce our basis. 

\medskip

\section{Case of singularity $1$}

The result of this section is published in \cite{Vi}.  Let us fix an $1$-singular point $o=(x^0_{kj})\in V$ such that $x^0_{ki} - x^0_{kj}\in \mathbb Z$.  We put $z_1= x_{ki} - x_{kj}$. Let $W$ be a sufficiently small neighborhood of the orbit  $\Gamma(o) = \Gamma(x^0)$ that is invariant with respect to the group $\Gamma$ and with respect to $\tau\in G$, where $\tau\in G$ is defined by $\tau(z_1)=-z_1$ and $\tau(z_i)=z_i$, $i>1$. In this case we have
$$
\mathcal L:= ev_o \circ \frac{\partial}{\partial z_1}\cdot z_1. 
$$
 In this case our basis has the following form
\begin{equation}\label{eq diff operator basis}
	\begin{split}
		D^1_{\sigma}:= \mathcal L \circ (\sigma + \tau(\sigma) ),\quad
		D^2_{\sigma'}:= \mathcal L \circ \frac{(\sigma' - \tau(\sigma') )}{z_1}, \quad \sigma, \sigma' \in \Gamma.
	\end{split}
\end{equation}
Here $D^1_{\sigma}$ is the sum of two local distributions $\mathcal L \circ \sigma$ at the point $\sigma^{-1}(o)$ and $\mathcal L \circ \tau(\sigma)$ at the point $(\tau(\sigma))^{-1}(o)$. The same holds for $D^2_{\sigma}$. 
We have the following equalities 
\begin{equation}\label{eq relations n=1}
D^1_{\tau(\sigma)} = D^1_{\sigma} \quad \text{and} \quad D^2_{\tau(\sigma')} = - D^2_{\sigma'}.
\end{equation}
From Theorem \ref{teor main n>1 class} it follows that 
 the local distributions $D^1_{\sigma},D^2_{\sigma'}$ up to relations (\ref{eq relations n=1}) form a basis for a $\mathfrak {gl}_n(\mathbb C)$-module. This $\mathfrak {gl}_n(\mathbb C)$-module was constructed  in \cite{Futorny} and it was  called the universal $1$-singular Gelfand-Tsetlin module. We can reformulate Theorem \ref{teor main n>1 class} in this case as follows. 

\medskip

\t\label{teor main} {\sl 
	The vector space spanned by $(D^1_{\sigma}, D^2_{\sigma'})$ up to relations (\ref{eq relations n=1}), where $\sigma,\sigma' \in \Gamma$,  is a $\mathfrak{gl}_n(\mathbb C)$-module. The action is given by Formulas (\ref{eq formula n=1, sing}).
	
}

\medskip

So we reproved Theorem $4.11$ from \cite{Futorny}. The explicit correspondence between notation in \cite{Futorny} and our notations can be deduced from  
$$
\mathcal {D}T(v+z):= D^2_{\sigma'}(T(v)), \quad \mathcal {D}^{\bar v}(F) = \frac{\partial}{\partial z_1}\Big|_{o} (F),
$$
where $\sigma'(v) = v+z$ and $T(v)$ and $F$ are holomorphic  function on $W$.

\bigskip

\noindent{\it Elizaveta Vishnyakova}

\noindent {Universidade Federal de Minas Gerais, Brazil}

\noindent{\emph{E-mail address:} \verb"VishnyakovaE@googlemail.com"}

\end{document}